\documentclass[12pt]{article}

\usepackage[colorlinks=true, pdfstartview=FitV, linkcolor=blue, citecolor=blue, urlcolor=blue]{hyperref}

\usepackage{amssymb,amsmath,amscd,array}
\usepackage{times, verbatim}
\usepackage{graphicx}
\usepackage{graphics}
\usepackage{epic}
\input dynkin.sty
\DynkinSmall
\DeclareFontFamily{OT1}{rsfs}{}
\DeclareFontShape{OT1}{rsfs}{n}{it}{<-> rsfs10}{}
\DeclareMathAlphabet{\mathscr}{OT1}{rsfs}{n}{it}
\let\mathcal=\mathscr

\usepackage{amssymb,amsmath, amscd}
\usepackage{times, verbatim}
\usepackage{graphicx}

\DeclareFontFamily{OT1}{rsfs}{}
\DeclareFontShape{OT1}{rsfs}{n}{it}{<-> rsfs10}{}
\DeclareMathAlphabet{\mathscr}{OT1}{rsfs}{n}{it}

\newtheorem{theorem}{Theorem}[section]

\newtheorem{con}[equation]{Conjecture}

\DeclareMathOperator{\End}{End}
\DeclareMathOperator{\Hom}{Hom}
\DeclareMathOperator{\Lie}{Lie}
\DeclareMathOperator{\Gal}{Gal}

\DeclareMathOperator{\SO}{SO}
\DeclareMathOperator{\Spin}{Spin}

\DeclareMathOperator{\SL}{SL}

\DeclareMathOperator{\disc}{disc}

\DeclareMathOperator{\Res}{Res}

\DeclareMathOperator{\Mass}{Mass}

\def\R{{\mathbb R}}
\def\F{{\mathbb F}}
\def\bC{{\mathbb C}}
\def\Q{{\mathbb Q}}
\def\Z{{\mathbb Z}}
\def\A{{\mathbb A}}

\title{Incoherent definite spaces and Shimura varieties}

\author{Benedict H Gross}

\begin{document}
\maketitle

\tableofcontents

\section{Introduction}

In \cite{G} we considered the infinite set of quaternion algebras $B(v)$ over a totally real field $k$ whose ramification locus has distance one from a finite set $\Sigma$ of places of $k$ which has odd cardinality. These algebras are indexed by the places $v$ of $k$; the algebra $B(v)$ is ramified at the set $\Sigma \cup \{v\}$ when $v$ is not contained in $\Sigma$, and is ramified at the set $\Sigma -\{v\}$ if $v$ is contained in $\Sigma$. When the set $\Sigma$ contains all the real places of $k$, each neighboring algebra $B(v)$ gives local information on a single Shimura curve $S$ which is defined over $k$. For example, if $v$ is a real place, the algebra $B(v)$ is split at $v$ and ramified at all other real places, and $\Sigma$-arithmetic subgroups of the multiplicative group of $B(v)$ give an analytic description of the points of $S$ over the quadratic extension $K_v = \bC$ of the completion $k_v = \R$ (See \cite{Sh} \cite{Sh2} for the construction of a canonical model for $S$ over $k \hookrightarrow \bC$ and \cite{DN} for the analytic description of $S$ at other real places). If $v = \frak p$ is a finite place where the curve has good reduction, the algebra $B(\frak p)$ is ramified at all real places and $\Sigma$-arithmetic subgroups of its multiplicative group give an analytic description of the points of $S$ which have "supersingular" reduction modulo $\frak p$ (cf. \cite[\S 11]{C}).

In this note, we generalize the notion of an odd set $\Sigma$ of places of $k$, containing all the real places (which can be viewed as an incoherent definite quaternion data over $k$) to incoherent definite orthogonal and Hermitian data over $k$. The former is a collection of local orthogonal spaces $V_v$ over the completions $k_v$ of fixed rank $n \geq 3$ and determinant $d$ in $k^*/k^{*2}$. We insist that the spaces $V_v$ are positive definite for all real places of $k$ (so the determinant $d$ is totally positive in $k^*$), that the Hasse-Witt invariants of $V_v$ are equal to +1 for almost all places $v$, and that the product of the Hasse-Witt invariants is equal to $-1$. There is no global orthogonal space over $k$ with these completions at all places, for then the product of the Hasse-Witt invariants would be equal to +1 by Hilbert's reciprocity law. However, for each place $v$ of $k$ there is a global orthogonal space $V(v)$ of rank $n$ and discriminant $d$ which is locally isomorphic to $V_u$ at all places $u \neq v$, and {\bf not} locally isomorphic to $V_v$ at $v$. The last condition determines the local orthogonal space when $v$ is finite; for a real place we also insist that the signature is equal to $(n-2,2)$. The global space $V(v)$ is determined up to isomorphism by its localizations, by the Hasse-Minkowski theorem. We call the global orthogonal spaces $V(v)$, indexed by the places $v$ of $k$, the neighbors of the incoherent definite data $\{V_v\}$.

An odd set $\Sigma$ of places of $k$ which contains all the real places gives incoherent definite orthogonal data of dimension $n = 3$ and determinant $d \equiv 1$. Indeed, for each place $v$, let $B_v$ be the local quaternion algebra over $k_v$ which is split if $v$ is not in $\Sigma$ and is ramified if $v$ is in $\Sigma$. Let $V_v$ be the local orthogonal space given by the norm form on elements of trace zero in $B_v$. The Hasse-Witt invariant of $V_v$ is the product of the Hasse invariant of the quaternion algebra $B_v$ with the class $(-1,-1)_v$ in the local Brauer group, so almost all of the Hasse-Witt invariants are equal to $+1$ and the product of Hasse-Witt invariants is equal to $-1$. The elements of trace zero in the nearby quaternion algebras $B(v)$ defined above give the neighboring global orthogonal spaces $V(v)$. 

There is a similar definition of an incoherent definite Hermitian data and its neighbors, where we fix the rank $n \geq 1$ of the space and a totally complex quadratic extension field $K$ of $k$. In this case, we only define a neighboring Hermitian space $V(v)$ at the places $v$ of $k$ which are not split in $K$. After working through the definitions, we show how incoherent definite orthogonal data of dimension $n \geq 3$ determines a Shimura variety $S$ of orthogonal type and dimension $n-2$, with field of definition $k$. We show how incoherent definite Hermitian data of dimension $n \geq 1$ determines a Shimura variety $S$ of unitary type and dimension $n-1$, with field of definition $K$. The special orthogonal groups of the neighbors of orthogonal data at real places $v$ can be used to describe $S$ over the complex quadratic extension $K_v$ of the completion $k_v$. Similarly, the unitary groups of the neighbors of Hermitian data at real places $v$ can be used to describe $S$ over the complex completions $K_v$. 

At primes $\frak p$ of $k$ where the orthogonal Shimura variety $S$ has good reduction, we propose to use the neighboring orthogonal space $V(\frak p)$ to describe the points over the unramified quadratic extension $K_{\frak p}$ of the completion $k_{\frak p}$ which reduce to a finite set of special points over the residue field $\F_{\frak p^2}$. When $k = \Q$, these special points should form the $0$-dimensional stratum of the supersingular locus. They should correspond to the moduli of orthogonal motives of rank $n$ and weight $2$ over $\F_{p^2}$ whose crystalline cohomology is a non-degenerate lattice of rank $n$ over $\Z_{p^2}$ with a semi-linear endomorphism $\phi$. The $\Z_p$ sublattice where $\phi = p$ should have rank $n$ and discriminant lattice isomorphic to $\F_{p^2}$ with its norm form.  

When the incoherent orthogonal data has dimension $n = 21$ over $\Q$ and is split at all finite primes, the Shimura variety $S$ parametrizes polarized $K3$ surfaces, and the stratification of the supersingular locus modulo $\frak p$ was first studied by Michael Artin \cite{A}. In that case, the $0$-dimensional stratum of special points consists of the moduli of supersingular $K3$ surfaces $X$ with Artin invariant $\sigma(X) = 1$. When $k \neq \Q$ we do not have a good definition of the special locus modulo $\frak p$, as there is no modular description of the points of $S$ modulo $\frak p$. But our special locus agrees with the superspecial locus, defined group-theoretically by Viehmann and Wedhorn \cite{ViW}.

Similarly at the primes $\frak p$ of $k$ which are inert in $K$ and where the unitary Shimura variety $S$ has good reduction, we propose to use the neighboring Hermitian spaces $V(\frak p)$ to describe the points over the completion $K_{\frak p}$ which reduce to a finite set of points in the special locus modulo $\frak p$. In the Hermitian case, the points of a Shimura variety closely related to $S$ have a modular interpretation, classifying abelian varieties with extra structure \cite{RSZ}, and the special points form the $0$-dimensional stratum of the supersingular locus \cite{VW}. A supersingular abelian variety $A$ is isogenous to a product of supersingular elliptic curves, and its modulus lies the the special locus precisely when $\dim \Hom(\alpha_p,A[p]) = \dim A$. This stratification of the supersingular locus for the moduli space of polarized abelian varieties  was first studied by Franz Oort (cf. \cite{O2} and the historical remarks in \cite{LO}). 

This paper represents my attempt to approach some of the work of Kudla and Rapoport on incoherent Eisenstein series and the arithmetic intersection of cycles on Shimura varieties \cite {KR1} \cite {KR2}. It is an expanded version of a letter that I wrote to Deligne in $2009$. Defining these Shimura varieties via incoherent quadratic and Hermitian spaces, rather than their neighboring special orthogonal and unitary groups at a real place $v$ of $k$, was suggested by the arithmetic conjectures in \cite{GGP}.

\section{Local and global orthogonal spaces}

In this section, we review the classification of orthogonal spaces over number fields and their completions. The main results are due to Minkowski, Witt, and Hasse. Some good references are \cite{L} \cite{MH} \cite [Ch IV]{S}.

We begin with a review of some invariants of an orthogonal space over a general field k, of characteristic not equal to $2$. 
An orthogonal space $(V,q)$ is a finite dimensional vector space $V$ over $k$ together with a quadratic form $q: V \rightarrow k$. The bilinear form
$$\langle v,w \rangle = q(v+w) - q(v) - q(w)$$
is symmetric, and we will always assume that it is non-degenerate (so gives an isomorphism from $V$ to its dual). The first invariant of the space $(V,q)$ is the rank, defined as the dimension of $V$. 

We can choose a basis $\{e_1.e_2, \ldots e_n\}$ of $V$ whose elements are pairwise orthogonal so that
$$q(\sum x_ie_i) = \sum a_i x_i^2.$$
The coefficients $a_i$ are all non-zero elements in $k$ and their product
$$d = a_1a_2 \cdots a_n$$
is well-defined in $k^*/k^{*2}$, independent of the orthogonal basis chosen. The class of $d$ in $k^*/k^{*2}$ is called the determinant of the quadratic space $(V,q)$, and is the second invariant we will use. If the rank is equal to $1$, the determinant $d$ determines the orthogonal space up to isomorphism -- the space is isometric to $k$ with the quadratic form $q(x) = dx^2$.

For two classes $a_i$ and $a_j$ in $ H^1(k, \mu_2) = k^*/k^{*2}$, we let $(a_i,a_j)$ denote their cup product in a $H^2(k, (\mu_2)^{\otimes 2}) = H^2(k, \mu_2) = Br_2(k)$. Here $Br(k) = H^2(k,\mathbb G_m)$ is the Brauer group of the field $k$. The Witt invariant of the quadratic space $(V,q)$ is defined by choosing an orthogonal basis and taking the sum in the Brauer group
$$w = \sum_{i < j} (a_i,a_j).$$
An important theorem, due to Witt, states that this class is again independent of the orthogonal basis chosen. This is the third invariant of $(V,q)$ we will need. 

If the rank of $(V,q)$ is equal to $2$, the determinant $d$ and the Witt invariant determine the orthogonal space up to isomorphism. In this case, the two dimensional $k$ vector space $V$ has the structure of a free module of rank one over the \'etale algebra $K = k[t]/(t^2 + d)$ (in two different ways). Indeed, any element in the special orthogonal group of $V$ which is not a scalar satisfies a quadratic polynomial whose roots generate the quadratic extension $K$. Picking one of the roots gives the action of $K$ on $V$ (the other root gives the conjugate action). Then $V$ is a free module of rank $1$ over $K$ with basis $e$, $q(e) = a$ is non-zero in $k$, and we have $q(\alpha.e) = a.N(\alpha)$ for all $\alpha$ in $K^*$, where $N:K^* \rightarrow k^*$ is the norm. The Witt invariant is determined by the class of $a$ in $Br_2(K/k) = k^*/N(K^*)$. Using the basis $1$ and $t$ for $K$ over $k$, and writing $\alpha = x + yt$ we see that the quadratic form on $V$ is given by $ax^2 + ady^2$ and that $V$ has Witt invariant $w = (a,ad) = (a,-d)$.

We call the choice of a homomorphism $i: K \rightarrow \End(V)$ of $k$-algebras such that $q(i(\alpha).v) = N(\alpha) q(v)$ for all $\alpha \in K^*$ and $v \in V$ an {\bf orientation} of the two dimensional space $V$. The orientations on $V$ form a principal homogenous space for the group $O(V)/\SO(V) = \mu_2$, and the group preserving the quadratic form on $V$ and a fixed orientation is  $\SO(V)$. When $K$ is a field, we can also give an orientation of $V$ by choosing one of the two isomorphisms
$(\Res_{K/k} \mathbb G_m)_{N = 1} \cong \SO(V)$ of non-split one dimensional tori over $k$. Indeed, the $K^*$-action on $V$ is determined by the action of the elements of norm $1$, so by the elements in $\SO(V)$. Finally, the choice of an orientation determines, and is determined by the choice of an isotropic line in the split orthogonal space $V \otimes K$ (where there are two isotropic lines). From an orientation and an element $\alpha$ in $\SO(V) \cong K^*_{N=1}$, we obtain an isotropic vector $e \otimes \alpha + \alpha(e) \otimes 1$ from any non-zero vector $e \in V$. Conversely, the choice of an isotropic line in $V \otimes K$ gives an isometric action of the group $K^*$ on this space, where $\alpha \in K^*$ acts by multiplication by $\alpha$ on the chosen isotropic line and by multiplication by $\alpha^{-1}$ on the other isotropic line. When $\alpha$ has norm $1$, this action preserves the space $V$ and gives an isomorphism $(\Res_{K/k} \mathbb G_m)_{N = 1} \cong \SO(V)$.

As a final example in in low rank, let $B = k + ki + kj + k(ij)$ be a quaternion algebra over $k$, with products $i^2 = a, j^2 = b, ij = -ji$. Let $V$ be the orthogonal space of dimension three over $k$, given by the norm form on elements of trace zero in $B$. Then
$$q(xi + yj + z(ij)) = -ax^2 -by^2 +abz^2.$$
The determinant of the $d(V,q)$ is $\equiv 1$ in $k^*/k^{*2}$, and the Witt invariant of $(V,q)$ is given by the sum $(a,b) + (-1,-1)$ in $Br_2(k)$.

When $k = k_v$ is a local field which is not isomorphic to $\bC$, the group $Br_2(k)$ is isomorphic to $\langle \pm1 \rangle$ and the cup-product of two classes in $k^*/k^{*2}$ is given by the Hilbert symbol $(a_i,a_j)_v = \pm1$. In this case, the Witt invariant is often referred to as the Hasse-Witt invariant, and is written as a product of signs:
$$\epsilon_v(V,q) = \prod_{i<j} (a_i,a_j)_v = \pm1.$$
When $k_v = \bC$ we put $\epsilon_v = +1$.

When $k_v = \R$, there is a further invariant of $(V,q)$ -- the signature $(r,s)$. Here $r$ is the dimension of a maximal positive definite subspace of $V$ and $s$ is the dimension of a maximal negative definite subspace of $V$. We have $r + s = \dim(V)$ and there is an orthogonal basis with
$$q(\sum x_ie_i) = x_1^2 + x_2^2 + \ldots + x_r^2 - x_{r+1}^2 - x_{r+2}^2 + \ldots - x_n^2.$$
We have
$$d_v \equiv (-1)^s ~~~~~ \epsilon_v = (-1)^{s(s-1)/2}.$$

Over the complex numbers, the dimension is a complete invariant of $(V,q)$ and over the real numbers, the signature is a complete isomorphism invariant of $(V,q)$. Over a non-Archimedean field, the dimension, determinant, and Hasse-Witt invariant are complete isomorphism invariants of $(V,q)$. The latter two invariants may be chosen arbitrarily in $k^*/k^{*2}$ and $\langle \pm 1 \rangle$, provided that $\dim(V) \geq 3$. In particular, for a fixed dimension  $n \geq 3$ and determinant $d$ there are precisely two isomorphism classes of orthogonal spaces over a non-Archimedean field, which are distinguished by the value of their Hasse-Witt invariant.

With these local invariants in hand, we can now describe the isomorphism classes of global orthogonal spaces $(V,q)$ over a number field $k$, with fixed dimension $n \geq 3$ and determinant $d$ in $k^*/k^{*2}$. At each real place $v$ choose a signature $(r_v,s_v)$ with $d_v \equiv (-1)^s$, and put $\epsilon_v = (-1)^{s(s-1)/2}$. At each finite place $v$, choose $\epsilon_v$ in $\langle \pm 1 \rangle$, with $\epsilon_v = +1$ for almost all $v$. The choice of local invariants gives us, for each place $v$, a orthogonal space over $k_v$ which is unique up to isomorphism. A necessary and sufficient condition for the existence of a global orthogonal space $(V,q)$ over $k$ of dimension of $n$ and determinant $d$, with these localizations is given by Hilbert's reciprocity law: $\prod_v \epsilon_v = +1$. In this case, the theorem of Hasse and Minkowski states that the global space $(V,q)$ is determined up to isomorphism over $k$ by its localizations. 

We may summarize these results as follows.

\begin{theorem}

Let $k$ be a number field, let $n \geq 3$ be an integer and let $d$ be a class in $k^*/k^{*2}$. For each place $v$ of $k$, let $(V_v,q_v)$ be a local orthogonal space over $k_v$ of dimension $n$ and determinant $d$. Assume that the local Hasse-Witt invariants $\epsilon_v = \epsilon_v(V_v,q_v)$ are equal to +1 for almost all places $v$. 

Then a necessary and sufficient condition for the existence of a global orthogonal space $(V,q)$ over $k$ with these localizations is 
$$\prod_v \epsilon_v = +1.$$
If a global orthogonal space $(V,q)$ exists with these localizations, it is unique up to isomorphism.

\end{theorem}

\section{Local and global Hermitian spaces}

Let $k$ be a number field, and let $K$ be a quadratic field extension of $k$. In this section, we review the classification of Hermitian spaces over $k$ and over the completions $k_v$. These results are due to Landherr, and a reference is \cite{MH}. We begin with the theory over a general field $k$ and a separable quadratic field extension $K$. Let $\tau$ be the non-trivial involution of $K$ over $k$.

A Hermitian space $(V,\phi)$ is a finite dimensional vector space over $K$ together with a Hermitian symmetric form $\phi: V \times V \rightarrow K$. We assume that $\phi$ is $K$-linear in the first variable and satisfies $\phi(w,v) = \phi(v,w)^{\tau}$. Then $\phi$ is $K$-anti-linear in the second variable and the quadratic form $q(v) = \phi(v,v)$ takes values in $k$. We also assume that $\phi$ is non-degenerate, so gives an isomorphism from $V$ to its dual. If $A = (\phi(e_i,e_j))$ is the matrix of the values of $\phi$ on a basis for $V$, then $A$ is Hermitian symmetric and $d = \det(A)$ is an element of $k^*$. Let $N(K^*)$ be the subgroup of $k^*$ consisting of those elements which are norms from $K^*$. Then the value of $d$ in the quotient group $k^*/N(K^*)$ is independent of the basis chosen, and is called the Hermitian determinant of $(V,\phi)$. A one dimensional Hermitian space is determined up to isomorphism by its Hermitian determinant -- the space is isomorphic to $K$ with form $\phi(z,w) = d. z w^{\tau}$.

One can also define Hermitian spaces when $K$ is the \'etale quadratic algebra $k \times k$ with involution $(a,b)^{\tau} = (b,a)$. In this case, for each dimension $n \geq 1$ there is a unique non-degenerate space up to isomorphism. The free $K$-module $V$ of rank $n$ has the form to $V = W + W'$, where he $W$ is a vector space of dimension $n$ over $k$ and $W'$ is its dual space. The pairing $\phi$ is defined by
$$\phi(w + w', u + u') = (\langle w,u' \rangle, \langle u,w' \rangle)$$
 In this case, all elements of $k^*$ are norms and $d \equiv 1$. 

If $k = k_v$ is a non-Archimedean local field, the dimension and the Hermitian determinant are complete isomorphism invariants of $(V,\phi)$. When $K = K_v$ is a field, the group $k^*/N(K^*)$ has order two and there are precisely two Hermitian spaces of each dimension $n \geq 1$. Moreover, when $K_v$ is the unramified quadratic field extension of $k_v$, the subgroup of norms are just the elements of $k^*$ with even valuation. If $K_v$ is not a field, there is a unique Hermitian space of each dimension. If $k_v = \R$ and $K_v = \bC$, there is a further invariant of a Hermitian space, its signature $(r,s)$. Here $r$ is the dimension of a maximal positive definite subspace and $s$ is the dimension of a maximal negative definite subspace. We have $r+s = \dim(V)$ and $d \equiv (-1)^s$. The signature is a complete isomorphism invariant in the real case.

Finally, let $k$ be a number field and let $K$ be a quadratic field extension of $k$. We want to describe all Hermitian spaces of a fixed dimension $n$ over $K$, up to isomorphism. At each real place $v$ of $k$ which is complex in $K$, choose a signature $(r_v,s_v)$ with $r_v + s_v = n$ and let $\epsilon_v = (-1)^{s_v}$. At each finite place $v$ of $k$ which is not split in $K$, choose a class $\epsilon_v$ in $k_v^*/N(K_v^*) = \langle \pm1 \rangle$ with $\epsilon_v = +1$ for almost all $v$. This choice gives us local Hermitian spaces of dimension $n$ for all places $v$, and the condition that these are the localizations of a fixed global space over $k$ is given by: $\prod_v \epsilon_v = +1$. When a global space exists, it is uniquely determined up to isomorphism by its local invariants. We summarize this result as follows.

\begin{theorem}
Let $k$ be a number field and let $K$ be a quadratic field extension of $k$. Fix a dimension $n \geq 1$. For each place $v$ of $k$ let $(V_v,\phi_v)$ be a local Hermitian space over $K_v$ of dimension $n$. Assume that the Hermitian determinant $d_v$ is a norm from $K_v^*$ for almost all places $v$, or equivalently that the Hilbert symbol $(d_v, \disc_{K/k})_v = +1$ for almost all $v$.

Then a necessary and sufficient condition for the existence of a global Hermitian space $(V,\phi)$ over $K$ with these localizations is given by:
$$\prod_v (d_v, \disc_{K/k})_v = +1$$
If a global Hermitian space $(V,\phi)$ exists with these localizations, it is unique up to isomorphism.
\end{theorem}

\section{Incoherent definite orthogonal and Hermitian data}

We henceforth assume that $k$ is a totally real number field, and in the Hermitian case, that $K$ is a totally complex quadratic field extension of $k$. Let $k^*_+$ denote the elements in $k^*$ which are positive at all real places.

Incoherent definite orthogonal data consists of a dimension $n \geq 3$, a determinant $d$ in $k^*_+/k^{*2}$, and for each place $v$ of $k$ a local orthogonal space $V_v$ of dimension $n$ and determinant $d$ over $k_v$. We assume that $V_v$ is positive definite for all real places $v$, that the Hasse-Witt invariant $\epsilon_v = \epsilon_v(V_v)$ is equal to $+1$ for almost all $v$, and that $\prod_v \epsilon_v = -1$.

Incoherent definite Hermitian data consists of a dimension $n \geq 1$ and for each place $v$ of $k$ a local Hermitian space $V_v$ over $K_v = K \otimes k_v$ of dimension $n$. We assume that $V_v$ is positive definite for all real places $v$, that the Hermitian determinant $d_v$ of $V_v$ is a norm from $K_v^*$ for almost all $v$, and that the number of places where $d_v$ is not a norm is odd. Equivalently, we assume that the Hilbert symbol $(d_v, \disc_{K/k})_v$ is equal to $+1 $ for almost all $v$ and that $\prod_v (d_v,\disc_{K/k})_v = -1$.

We will use the notation $\{V_v\}$ for incoherent definite data, either orthogonal or Hermitian.

Given incoherent definite orthogonal data over $k$, we construct for each place $v$ of $k$ a global orthogonal space $V(v)$ of dimension $n$ and determinant $d$ over $k$ which we call the neighbor at $v$. This space has localizations $V_u$ for all places $u \neq v$. At the place $v$ we insist that
$$\epsilon_v(V(v)) = - \epsilon_v(V_v)$$
This determines the localization of $V(v)$ at a finite place $v$. At a real place, we insist further that the signature of $V(v)$ at $v$ is equal to $(n-2,2)$. Since the product of local epsilon factors is now equal to $+1$, there is a global space $V(v)$ with these localizations, which is unique up to isomorphism. Finally, for each place $v$ of $k$ we let $G(v) = \SO(V(v))$ be the special orthogonal group of the neighboring orthogonal space.

Given incoherent definite Hermitian data over $k$ relative to the quadratic extension $K$, we construct for each place $v$ of $k$ which is not split in $K$ a global Hermitian space $V(v)$ of dimension $n$ over $K$ which we call the neighbor at $v$. This space has localizations $V_u$ for all places $u \neq v$. At the place $v$ we insist that
$$(d(V(v)), \disc_{K/k})_v = - (d(V_v), \disc_{K/k})_v$$
This determines the localization of $V(v)$ at a finite place $v$. At a real place, we insist further that the signature of $V(v)$ at the place $v$ is equal to $(n-1,1)$. Since the product of local Hilbert symbols is now equal to $+1$, there is a global Hermitian space $V(v)$ with these localizations, which is unique up to isomorphism. Finally, for each place $v$ of $k$ which is not split in $K$, we let $G(v) = U(V(v))$ be the unitary group of the neighboring Hermitian space.

 \section{Incoherent data in codimension one}

Let  $k$ be a totally real field and let $ \{V_v\}$ be incoherent definite orthogonal data of dimension $n \geq 4$ and determinant $d$, and let $a$ be any class in $k^*_+/k^{*2}$. Then $a$ is represented by all of the local quadratic spaces $V_v$. Let $U_v$ denote the orthogonal complement of the corresponding one dimensional subspace of $V_v$, so that
$$V_v = U_v \oplus \langle a \rangle.$$
Then $\{U_v\}$ is incoherent definite orthogonal data of dimension $n-1$ and determinant $ad$. Clearly the subspace $U_v$ is positive definite at all real places $v$. But we also have the formula for the Hasse-Witt invariants
$$\epsilon_v(U_v) = \epsilon_v(V_v) \cdot (da,a)_v.$$
with $(da,a)_v = +1$ for almost all $v$ and $\prod_v (da,a)v = +1$. Hence $\epsilon_v(U_v)= +1$ for almost all $v$ and $\prod_v \epsilon_v(U_v) = -1.$ It follows from a calculation of Hasse-Witt invariants that the neighboring spaces $U(v)$ all have the property that 
$$U(v) \oplus \langle a \rangle \simeq V(v).$$

Conversely, if we are given incoherent definite data $\{U_v \subset V_v\}$  of codimension $1$, with determinants $d(U)$ and $d(V)$, let $a = d(V)/d(U)$ in $k^*_+/k^{*2}$. Then $U_v^{\perp} = \langle a \rangle$ is coherent definite orthogonal data of dimension $1$. Hence we have shown that

\begin{theorem}
For incoherent definite orthogonal data $V = \{V_v\}$ of dimension $n \geq 4$, the classes of incoherent definite orthogonal data $U = \{U_v\}$ of codimension one in $V$ correspond bijectively to classes $a$ in $k^*_+/k^{*2}$.
\end{theorem}

There is a similar result when the dimension of $V$ is $3$, but one has to restrict to classes $a$ in $k^*_+/k^{*2}$ which are locally representable at the finite set of places $v$ where the space $V_v$ is anisotropic. Also, incoherent definite orthogonal data $\{U_v\}$ of dimension $2$ and determinant $d$ is better viewed as incoherent definite Hermitian data of dimension $1$ over the imaginary quadratic extension $K = k(\sqrt{-d})$ of $k$.
In the Hermitian case, we have the following.

\begin{theorem}
Let $K$ be an imaginary quadratic extension of the totally real field $k$. For incoherent definite Hermitian data $\{V_v\}$ of dimension $n \geq 2$ over $K$, the classes of incoherent definite Hermitian data $\{U_v\}$ of codimension one in $V$ correspond bijectively to classes $a$ in $k^*_+/N(K^*)$.
\end{theorem}

Indeed, the subspaces $U_v^{\perp} = \langle a_v \rangle$ give coherent definite Hermitian data of dimension $1$. Since $\prod_v (a_v, \disc{K/k}) = +1$, the class of $(a_v)$ is trivial in the quotient group $\A_k^*/N(\A_K^*)k^*$, which has order $2$. Hence the local classes $a_v$ can be written as a product of a global class $a$ and a norm from $K_v^*$, and $a$ is well-defined in the quotient group $k^*_+/N(K^*)$.

\section{Orthogonal lattices}

In this section, $k$ is a local non-Archimedian field, with ring of integers $A$, uniformizing parameter $\pi$, and finite residue field $\F = A/\pi A$ of order $q$. We will assume further that the characteristic of $\F$ is odd.

Let $V = (V,q)$ be a non-degenerate orthogonal space over $k$ of dimension $n$, with associated bilinear form $\langle v,w \rangle = q(v+w) - q(v) - q(w)$. A lattice $L \subset V$ is by definition a free $A$-module of rank $n$ which spans $V$ over $k$, such that the bilinear form $\langle , \rangle$ takes integral values on $L$. The dual lattice $L^{\vee}$ is the free $A$-module of vectors in $V$ with integral pairing against all elements of $L$. Then $L$ is a submodule of $L^{\vee}$ and the quotient $L^{\vee}/L$ has finite length, by the non-vanishing of the determinant $d$ of $V$. 

We say that the lattice $L$ is {\bf non-degenerate} (some say self-dual) if $L = L^{\vee}$. A necessary condition for the existence of a non-degenerate lattice in $V$ is that the determinant of $V$ has even valuation, so $d$ lies in the subgroup $A^*k^{*2}/k^{*2}$ of $k^*/k^{*2}$. In this case, we say that the determinant of $V$ is equivalent to a unit. This condition is clearly sufficient for the existence of a non-degenerate lattice in $V$ when the dimension $n = 1$, but is not sufficient when $n \geq 2$. 

To see this in the two dimensional case, let $D$ be a unit in $A^*$ which is not a square, so $K = k(\sqrt D)$ is the unramified quadratic extension of $k$. Let $V = K$ with quadratic form $q(v) =  N(v)$. Then $q(x,y) = x^2 - Dy^2$, the determinant of $V$ is equivalent to the unit $-D$ in $k^*/k^{*2}$ and $\epsilon(V) = +1$. In this case $V$ contains the non-degenerate lattice $A_K = A + \sqrt D A$ given by the ring of integers in $K$. On the other hand, if $V' = K$ with quadratic form $q(v) = \pi.N(v)$, then the determinant of $V'$ is still equivalent to the unit $-D$ but now $\epsilon(V') = -1$. In this case, there are no non-degenerate lattices in $V'$, but $L = A_K$ is a lattice with dual lattice $L^{\vee} = (\pi)^{-1}A_K$. The quadratic space $L^{\vee}/L$ has rank $2$ over the residue field $\F$ and is isomorphic to the quadratic extension field $A_K/\pi A_K$ of order $q^2$ with its norm form. More generally, we have the following.

\begin{theorem}
Let $V$ be an orthogonal space over $k$ of dimension $n$ whose determinant $d$ is equivalent to a unit.

If $\epsilon(V) = +1$ then $V$ contains a non-degenerate lattice $L$. All such lattices are conjugate under the action of the group $\SO(V)$.

If $\epsilon(V) = -1$ then there are no non-degenerate lattices in $V$. However, $V$ contains lattices $L$ with the property that $L^{\vee}/L$ is a two dimensional orthogonal space over the residue field $\F$, isomorphic to the quadratic extension field of $\F$ with its norm form. All such lattices are conjugate under the action of the group $\SO(V)$.
 
 \end{theorem}
 
 This is well-known -- the lattices $L$ in the theorem are called maximal lattices by Eichler, as they are maximal with the property that the bilinear form takes values in $A$. We will sketch a proof which involves embedding two dimensional spaces $W$ in $V$, which will be useful later in this paper. The theorem is easily proved (using the construction above) when the dimension $n$ of $V$ is $\leq 2$, so we will assume that $n \geq 3$. 
 
 \begin{theorem}
Assume that the dimension of $V$ satisfies $n \geq 3$ and that the determinant $d$ of $V$ is equivalent to a unit. Let $D$ be a unit which is not a square, let $K$ be the unramified quadratic extension $k(\sqrt D)$ and let $W$ be an orthogonal space of dimension $2$ and determinant $-D$, which satisfies
$$\epsilon(W) = \epsilon(V).$$
Then $W$ embeds isometrically as a subspace of $V$, and all embeddings are conjugate under $\SO(V)$. The orthogonal complement $U$ of $W$ in $V$ has invariants
$$\dim(U) = n-2~~~~~d(U) \equiv -dD~~~~~~\epsilon(U) = +1$$
\end{theorem}

Note that the isomorphism class of the orthogonal complement $U$ depends only on the dimension and discriminant of $V$. This is similar to what happens when $k = \R$, for the signatures we specified. The space $V$ has signature $(n,0)$ or $(n-2,2)$, so has $d \equiv 1$. If we take $K = \bC$, so $D \equiv -1$, and insist that the two dimensional space $W$ has determinant $-D$ and the same Hasse-Witt invariant as $V$, then $W$ has signature $(2,0)$ in the first case and signature $(0,2)$ in the second case. Hence $W$ embeds in $V$, and in both cases the orthogonal complement $U$ has signature $(n-2,0)$.

To prove Theorem $6.2$, we will use some standard results \cite[Ch IV]{S} on which elements of $k^*$ are represented by quadratic forms. We begin with the case when the dimension of $V$ is equal to $3$. If a two dimensional space $W$ of determinant $-D$ embeds into $V$, its orthogonal complement $U$ would be one dimensional, of determinant $-dD$. But the quadratic form $q$ on $V$ represents the unit $-dD$, as it represents all elements of $k^*$ except perhaps $-d$, and the discriminant $D$ of $K$ is not equivalent to $1$. Taking the orthogonal complement of this line spanned by a vector $v$ with $q(v) = -dD$ gives a two dimensional subspace $W$ of $V$ which has determinant $-D$. The quadratic form on $W$ is given by $ax^2 - aDy^2$ and on $V$ is given by $ax^2 -aDy^2 - Ddz^2$.  Since $-D$ and $-dD$ are both units and the residue characteristic is odd, $(-D,-dD) = +1$ and the Hasse-Witt invariants of both $W$ and $V$ are both equal to $(a,-D)$,
This proves our claim on the embedding when the dimension of $V$ is three. The orthogonal group of $V$ acts transitively on these two dimensional subspaces, by Witt's extension theorem. Since the stabilizer contains a reflection, the special orthogonal group also acts transitively. 

One can also prove the special case when the dimension of $V$ is $3$ and the determinant $d \equiv 1$ using the theory of local quaternion algebras. Indeed, in this case the space $V$ is given by the elements of trace zero in a quaternion algebra $B$ over $k$, and the quadratic form is the norm form. The unramified quadratic extension $K$ of $k$ embeds as a subfield of $B$, and its orthogonal complement $W$ is the desired two dimensional subspace of $V$. We obtain an orientation of the plane $W$ in $V$ by choosing one of the two vectors on the line $W^{\perp}$ with $q(v) = -D$, which also gives $W$ a $K = k.1 + k.v$  linear structure.

Next consider the case when the dimension of $V$ is equal to four. The quadratic space $V$ represents every non-zero class $a$ and the three dimensional orthogonal complement of the line spanned by a vector $v$ with $q(v) = a$ has determinant $ad$. To obtain an embedding of $W$, we need to check that this three dimensional space represents $-aD$. But any three dimensional quadratic space of determinant $ad$ represents all classes not equivalent to $-ad$,  Hence we have an embedding of both orthogonal spaces of dimension $2$ and determinant $-D$, except in the case when the determinant of $V$ is equivalent to $D$.
In that case, we can construct a four dimensional space $V$ of determinant $D$ from a two dimensional space $W$ of determinant $-D$ by taking the direct sum of $W$ with a hyperbolic plane. The quadratic form on $V$ has the form $ax^2  - aDy^2 + z^2 - t^2$, so the Hasse-Witt invariant of both $W$ and $V$ are both equal to $(a,-D)$.

When the dimension of $V$ is at least five, then the quadratic form on $V$ represents any non-zero class $a$ and the orthogonal complement of the line spanned by a vector $v$ with $q(v) = a$ represents $-aD$. Hence we may embed both quadratic spaces $W$ of dimension $2$ and determinant $-D$ in $V$. The special orthogonal group of $V$ acts transitively on the subspaces $W$ with $\epsilon(W) = \epsilon(V)$, by Witt's extension theorem. Hence the orthogonal complement $U$ of $W$ is well-defined up to isomorphism over $k$. Since $V = W \oplus U$ we find 
$$d(U) \equiv -dD ~~~~ \epsilon(U) = +1$$
as claimed.

We can now prove Theorem $6.1$ by an induction on the dimension of $V$. If $\epsilon(V) = +1$, we embed the two dimensional subspace $W$ with determinant $-D$ and $\epsilon(W) = +1$.
We have seen that there is a non-degenerate lattice in $W$, and by induction there is a non-degenerate lattice in its orthogonal complement $U$. The direct sum of these two lattices gives a non-degenerate lattice in $V$. If $\epsilon(V) = -1$ we embed the two dimensional subspace $W$ with determinant $-D$ and $\epsilon(W) = -1$. By induction there is a non-degenerate lattice in the orthogonal complement $U$, and we have seen that there is a lattice $L$ in $W$ with $L^{\vee}/L$ isomorphic to the quadratic extension of $\F$ with its norm form. The direct sum of these lattices gives the desired lattice in $V$.

The fact that the group $\SO(V)$ acts transitively on the lattices described in Theorem $5.1$  is proved in \cite{Sh} and \cite{GHY} They are the maximal lattices in $V$, where the bilinear form is integral. In the case when  $\epsilon = +1$ and the lattice $L$ is non-degenerate, the stabilizer of $L$ is a hyperspecial maximal compact subgroup of $\SO(V)$. When $n = 2m+1$ is odd, its reductive quotient is the split group $\SO_{2m+1}$ over $\F$. When $n = 2m$ is even and $(-1)^md \equiv 1$  its reductive quotient is the split group $\SO_{2m}$ over $\F$. When $n = 2m$ is even and $(-1)^md \equiv D$ its reductive quotient is the quasi-split group $\SO_{2m}^{\epsilon}$ over $\F$.

In the case when $\epsilon = -1$, let $T = \SO_2^\epsilon$ be the one dimensional torus over $\F$ which is split by the quadratic extension of $\F$. Then the stabilizer of $L$ is a parahoric subgroup of $\SO(V)$ whose reduction is as follows. When $n = 2m+1$ is odd, the reductive quotient is the group $T \times \SO_{2m-1}$ over $\F$. When $n = 2m$ is even and $(-1)^md \equiv 1$, the reductive quotient is the group $T \times \SO_{2m-2}^{\epsilon}$ over $\F$. When $n = 2m$ is even and $(-1)^md \equiv D$, the reductive quotient is the group $T \times \SO_{2m-2}$ over $\F$. When $n \geq 3$ the parahoric subgroup $\SO(L)$ is not hyperspecial. It is maximal except when $2m = 4$ and $d \equiv D$, when it is an Iwahori subgoup of $\SO(V) \cong (\Res_{K/k}\SL_2)/\mu_2$. See \cite{GHY} for details.

\section{Hermitian lattices}

In this section, we again assume that $k$ is a local non-Archimedian field, with ring of integers $A$, uniformizing parameter $\pi$, and finite residue field $\F = A/\pi A$ of order $q$. Let $K = k(\sqrt D)$ be the unramified quadratic extension of $k$, where $D$ is a unit in $A^*$ which is not a square, and let $\tau$ be the non-trivial involution of $K$ which fixes $k$. We make no assumption on the characteristic of $\F$.

Let $V$ be a non-degenerate Hermitian space over $K$ of dimension $n$, with associated Hermitian symmetric form $\phi (v,w)$. A lattice $L \subset V$ is by definition a free $A_K$-module of rank $n$ which spans $V$ over $K$, such that the Hermitian form $\phi$ takes integral values on $L$. The dual lattice $L^{\vee}$ is the free $A_K$-module of vectors in $V$ with integral pairing against all elements of $L$. Then $L$ is an $A_K$-submodule of $L^{\vee}$ and the quotient $L^{\vee}/L$ has finite length, by the non-vanishing of the Hermitian determinant $d$ of $V$. We say that the lattice $L$ is {\bf non-degenerate} (some say self-dual) if $L = L^{\vee}$. A necessary condition for the existence of a non-degenerate lattice in $V$ is that the determinant has even valuation, so $d$ lies in the subgroup $N(K^*)$ of index $2$ of $k^*$. In this case, we say that the determinant of $V$ is a norm. 

When the dimension of $V$ over $K$ is equal to $1$ and the determinant is not a norm, then the Hermitian space is isomorphic to $K$ with form $\phi(x,y) = \pi xy^{\tau}$. In this case, there are no non-degenerate lattices in $V$. However the lattice $L = A_K$ has dual lattice $L^{\vee} = \pi^{-1}A_K$, and the quotient module $L^{\vee}/L$ is isomorphic to the quadratic extension field of $\F$ with its standard Hermitian form.

\begin{theorem}

Let $V$ be a non-degenerate Hermitian space of dimension $n \geq 1$ over the unramified quadratic extension $K = k(\sqrt D)$ of $k$.

If the Hermitian determinant $d$ of $V$ is a norm, then $V$ contains a non-degenerate lattice $L$. All such lattices are conjugate under the action of the unitary group $U(V)$.

If the Hermitian determinant $d$ of $V$ is not a norm, then there are no non-degenerate lattices in $V$. However, $V$ contains lattices $L$ with the property that $L^{\vee}/L$ is a one dimensional Hermitian space over the quadratic extension of the residue field $\F$. All such lattices are conjugate under the action of the unitary group $U(V)$.

\end{theorem}

The proof is similar to the orthogonal case, but simpler as it involves the embedding of a suitable one dimensional subspace $W$ in $V$. Let $W$ be a Hermitian space of dimension $1$ with the property that
$$d(W) \equiv d(V)$$
in $k^*/N(K^*)$. It is then a simple matter to show that $W$ embeds isometrically as a subspace of $V$, and that all such subspaces are conjugate under the action of $U(V)$. The orthogonal complement of $W$ in $V$ has Hermitian determinant which is a norm, so by induction contains non-degenerate lattice. The direct sum of that lattice with the lattice constructed above (in dimension $1$) gives the desired lattice in $V$. 

The fact that the unitary group of $V$ acts transitively on lattices of this type follows from the theory of Iwahori and Bruhat-Tits on parahoric subgroups. Indeed, when the determinant of $V$ is a norm, the stabilizer of a non-degenerate lattice $L$ is a hyperspecial maximal compact subgroup of the locally compact group $U(V)$ whose reductive quotient is the group $U_n$ over $\F$. When the determinant of $V$ is not a norm, the stabilizer is again a maximal compact subgroup, but is not hyperspecial once $n \geq 2$. Its reductive quotient is the group $U_1 \times U_{n-1} = T \times U_{n-1}$ over $\F$.

\section{The homogeneous spaces $X$, $Y$, and $Z$}

We continue with the notation of the previous two sections: $k$ is a non-Archimedian local field with ring of integers $A$ and finite residue field $\F = A/\pi A$. In the orthogonal case, we assume that the characteristic of $\F$ is odd. Let $K = k(\sqrt D)$ be the unramified quadratic extension of $k$, and let $A_K$ be the ring of integers of $K$. 

Our aim is to generalize the following construction, where $k = \R$ and $K = \bC$. If $V$ is a real orthogonal space of dimension $n \geq 3$ and signature $(n-2,2)$, the set $X = \mathscr D^{\pm}$ of oriented negative definite two planes $W$ in $V$ is a complex manifold of dimension $n-2$ which admits a transitive action of the group $G =\SO(V)$. Indeed, an oriented negative definite two plane $W$ determines a homomorphism
$h: T \ \rightarrow \SO(V)$ over $\R$, where $T$ is the one dimensional torus $(\Res_{\bC/\R} \mathbb G_m)_{N = 1}$, and the centralizer of $h$ is the compact subgroup
$$H = \SO(W) \times \SO(W^{\perp}).$$
The conjugacy class of $h$ is isomorphic to the quotient space $X = G/H$. This quotient is isomorphic to two copies of the symmetric space of the connected real group $\SO(n-2,2)^+$, as the subgroup $H$ has index two in its normalizer $S(O(W) \times O(W^{\perp}))$ (which is a maximal compact subgroup of $G$). Hence every coset for the normalizer breaks into two cosets for the subgroup $H$. The space $X$ has a Hermitian symmetric structure, as the tangent space to the identity coset of $G/H$ is the orthogonal representation $\Hom(W, W^{\perp}) \cong W \otimes W^{\perp}$ of $H$, which has a complex Hermitian structure from the chosen orientation of $W$. The involution of $X$ which reverses the orientation of $W$ switches the two components and is anti-holomorphic. Finally, the complex manifold $X$ embeds as an open $G$-orbit (with two components) in the complex quadric $Q(\bC)$ defined by the vanishing of the quadratic form on the projective space of $V \otimes \bC$.

In the non-Archimedean case, we start with an orthogonal space $V$ over $k$ of dimension $n \geq 3$, with unit determinant $d$ and Hasse-Witt invariant $\epsilon = -1$. Let $X$ be the set of oriented two planes $W$ in $V$ with determinant $-D$ and satisfying $\epsilon(W) = -1$. We have seen that such oriented two planes exist, and are permuted transitively by $G = \SO(V)$. An oriented two plane $W$ determines a homomorphism $h: T \ \rightarrow \SO(V)$ over $k$, where $T$ is the one dimensional torus $(\Res_{K/k} \mathbb G_m)_{N = 1}$, and the centralizer of $h$ is the closed subgroup
$$H = \SO(W) \times \SO(W^{\perp}).$$
The conjugacy class $X$ of such homomorphisms $h$, which can be identified with the coset space $G/H$ once an oriented plane $W$ has been chosen, has the structure of a homogeneous analytic manifold over $k$. Indeed, if $G$ is a $k$-analytic group and $H$ is a $k$-analytic closed subgroup, the coset space $G/H$ has the structure of a $k$-analytic manifold, satisfying the expected universal properties, with tangent space $\Lie(G)/\Lie(H)$ at the identity \cite[Ch III,\S 12]{S2}. In this case, the tangent space of $X$ at the oriented plane $W$ is the $k$-vector space $\Hom(W,W^{\perp}) = W^{\vee} \otimes W^{\perp} \cong W \otimes W^{\perp}$, which affords an orthogonal representation of the stabilizer. Since the subspace $W$ has $K$-linear structure given by the orientation, the tangent plane at each point is a $K$-vector space. Hence $X$ is a $K$-analytic manifold of dimension $n-2$.

 As in the real case, there is a fixed point free involution of $X$, which acts $K$-anti-linearly on the tangent spaces and commutes with the action of $\SO(V)$. This is given by the action of normalizer of $h$ in $G$, which contains the centralizer with index $2$ and maps each homomorphism $h : T \rightarrow \SO(V)$ to its composition with the inverse homomorphism on $T$. Finally, the analytic manifold $X$ embeds as an open $G$-orbit in the $K$-points of the smooth quadric $Q$ of dimension $n-2$  defined by the vanishing of the quadratic form $q$ on the projective space of $V \otimes K$. The embedding sends the oriented two plane $W = K.e$ in $V$ spanned by the orthogonal vectors $e$ and $f = \sqrt D.e$ to the line in $V \otimes K$ spanned by the vector $e \otimes 1 + f \otimes (\sqrt D)^{-1}$. Since $q(f) = N(\sqrt D)q(e) = -Dq(e)$ this line is isotropic, and defines a point on $Q(K)$.

Unlike the real case, the subgroup $H$ need not be compact. It is compact only when the factor $\SO(W^{\perp})$ is compact, which can only occur when the dimension of $V$ is three or four. To obtain a compact subgroup of $H$, we need to choose a non-degenerate lattice $L$ in the subspace $W^{\perp}$ (which has unit determinant and $\epsilon(W^{\perp}) = +1)$. Writing $W = K.e$ with $q(e) = \pi$ we obtain a canonical lattice $A_K.e$ in $W$ and the product
$$J = \SO(A_K.e) \times \SO(L)$$
is a hyperspecial maximal compact subgroup of $H$. The choice of $L$ also gives us a maximal lattice $\Lambda = A_K.e  + L$ in $V$, and the orientation of the two plane $W$ gives an orientation of the two dimensional orthogonal space $\Lambda^{\vee}/\Lambda$ over the residue field $\F = A/\pi A$.

We define $Y$ as the set of pairs $(W,L)$, where $W$ is an oriented two plane with determinant $-D$ and Hasse-Witt invariant $\epsilon(W) = -1$ and $L$ is a non-degenerate lattice in $W^{\perp}$. The group $G = \SO(V)$ acts transitively on the set $Y$, which can be identified with the coset space $G/J$, once a fixed pair $(W,L)$ with stabilizer $J$ has been chosen. There is s a fixed point free involution of $Y$ given by the normalizer of $J$ in $G$, which contains $J$ with index $2$. This changes the orientation of the plane $W$ and preserves the lattice $L$.

We define $Z$ as the set of maximal lattices $\Lambda$ in $V$, together with an orientation of the orthogonal space $\Lambda^{\vee}/\Lambda$ over $A/\pi A$. The group $G = \SO(V)$ acts transitively on the set $Z$, and the stabilizer of a fixed maximal lattice with orientation is a parahoric subgroup $N$ (which is maximal if the dimension of $V$ is at least $5$). The reductive quotient of $N$ is a group of the type $\SO(2) \times \SO(n-2)$ \cite{GHY}, and is isomorphic to the reductive quotient of the hyperspecial maximal compact subgroup $J$ of $H$. We can identify the set $Z$ with the homogeneous space $G/N$ once a fixed oriented maximal lattice $\Lambda$ with stabilizer $N$ has been chosen. There is a fixed point free involution of $Z$ given by the normalizer of $N$ in $G$, which contains $N$ with index $2$ and is a maximal compact subgroup once the dimension of $V$ is at least $5$. This preserves the lattice $\Lambda$ and changes the orientation of $\Lambda^{\vee}/\Lambda$.

The map that takes the pair $(W,L)$ to the oriented two plane $W$ is a surjective $G$-equivariant map
$$Y \rightarrow X \hookrightarrow Q(K)$$
which is usually of infinite degree. The map that takes the pair $(W,L)$ to the oriented maximal lattice $\Lambda = A_K.e + L$ is a $G$-equivariant map 
$$f: Y \rightarrow Z.$$
The map $f$ is a morphism in the category of $G$ sets which commutes with the involutions we have defined on $Y$ and $Z$. Since $G$ acts transitively on $Z$, it is surjective. We will see in the next section that the fibers of $f$, which are all isomorphic to the quotient space $N/J$, have a simple $A_K$ analytic structure. 

The parahoric subgroups of the simply-connected cover $\Spin(V)$ of $\SO(V)$ can be read off the affine diagram, together with its Frobenius automorphism. In the case when the dimension $n= 2m + 1 \geq 7$ is odd, the determinant is a unit, and $\epsilon(V) = -1$, we have the diagram 
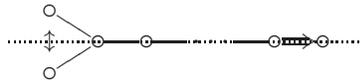
\[
\Dynkin{B*}{n[2]}{}
\]
The parahoric $N$ corresponds to the two vertices interchanged by Frobenius at the left end of the diagram.
In the case when the dimension $n = 2m \geq 8$ is even, the determinant is $\equiv (-1)^m$ and $\epsilon(V) = -1$ we have the diagram 
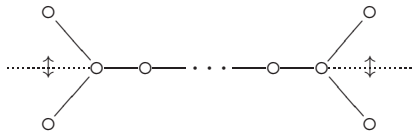
\[
\Dynkin{D*}{n[2]}{}
\]
and the parahoric $N$ corresponds to two vertices interchanged by Frobenius at the left end of the diagram.
In the case when the dimension $n = 2m$ is even and the determinant is $\equiv (-1)^mD$ we have the diagram 
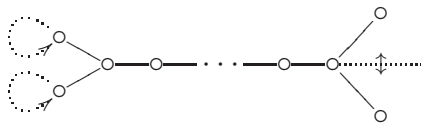
\[
\Dynkin{D*}{n[-2]}{}
\]
When $\epsilon(V) = -1$ the parahoric $N$ corresponds to the two vertices interchanged by Frobenius at the right end of the diagram. 

In all cases, when the dimension of $V$ is at least $5$, the reductive quotient of $N$ has center isomorphic to the non-split torus $T$ of dimension $1$ over $\F$. The finite group $T$ lifts to $N$, this lifting is well-defined up to conjugacy, and the centralizer of any lifting is conjugate to the subgroup $J$.

A similar theory works in the Hermitian case. First consider the situation over the real numbers. If $V$ is a complex Hermitian space of dimension $n$ and signature $(n-1,1)$, then the set $X$ of negative lines $W$ in $V$ is a complex manifold of dimension $n-1$ which admits a transitive action of the group $G = U(V)$. Indeed, such a line determines a homomorphism $h: T \rightarrow U(V)$ over $\R$, and the centralizer of $h$ is the compact subgroup
$$H = U(W) \times U(W^{\perp}).$$
The quotient $X \cong G/H$ is isomorphic to the symmetric space of $U(n-1,1)$. It has a Hermitian symmetric structure, as the tangent space to the identity is the representation $W \otimes W^{\perp}$ of $H$, of dimension $n-1$ over $\bC$. Finally, the complex manifold embeds as an open $G$ orbit in the projective space of $V$ over $\bC$.

In the non-Archimedian case, we start with a Hermitian space $V$ over $K$ of dimension $n \geq 1$ whose Hermitian determinant is not a norm. Let $X$ be the set of lines $W$ in $V$ with Hermitian determinant equal to the Hermitian determinant of $V$. We have seen that such lines exist and are permuted transitively by $G = U(V)$. Such a line determines a homomorphism $h: T \rightarrow U(V)$, where $T$ is the one dimensional torus $(\Res_{K/k} \mathbb G_m)_{N = 1}$, and the centralizer of $h$ is the closed subgroup
$$H = U(W) \times U(W^{\perp}).$$
The conjugacy class $X$ of such homomorphisms $h$ is isomorphic to the quotient space $G/H$ and $X$ embeds as an open $G$ orbit in the $K$ valued points of the projective space of $V$.

The subgroup $H$ can only be compact when the dimension of $V$ is one or two. To obtain a suitable compact subgroup of $H$, we need to choose a non-degenerate lattice $L$ in the subspace $W^{\perp}$. Writing $W = K.e$ with $\phi(e,e) = \pi$, we obtain a canonical lattice $A_K.e$ in $W = K.e$. The product
$$J = U(A_K.e) \times U(L)$$
is a hyperspecial maximal compact subgroup of $H$ and  
the homogeneous space $Y = G/J$ indexes the lines $W$ in $V$ with Hermitian determinant equal to the determinant of $V$ together with a non-degenerate lattice $L$ in $W^{\perp}$. The map
$$Y \rightarrow X \hookrightarrow \mathbb P(V\otimes K)$$
which assigns to a pair $(W,L)$ the line $W$ is usually of infinite degree.

Let $Z$ be the set of maximal lattices $\Lambda$ in $V$. The map 
$$f: Y \rightarrow Z$$
associates to the line $W$ with the correct Hermitian determinant and the non-degenerate lattice $L$ in $W^{\perp}$ the maximal lattice $ \Lambda = A_K.e + L$. The group $G = U(V)$ acts transitively on the space $Z$ of maximal lattices and the map $f$ is $G$-equivariant, hence surjective. Let $N = U(\Lambda)$ be the stabilizer of the maximal lattice $\Lambda$. This is a parahoric subgroup of $G$ whose reductive quotient is isomorphic to $U(1) \times U(n-1)$ over the residue field \cite {GHY} and is isomorphic to the reductive quotient of the compact subgroup $J$. We can identify the space $Z$ with the quotient space $G/N$, once a maximal lattice has been chosen. In the next section we will see that the fibers of the map $f$, which are all isomorphic to the quotient space $N/J$, have a simple $A_K$ analytic structure.

The parahoric subgroups of $U(V)$ can be read off the affine diagram of $SU(V)$ with its automorphism given by Frobenius. When the dimension $n = 2m \geq 4$ of $V$ is even and the determinant is not a norm, we have the diagram
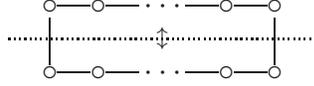
\[
\Dynkin{A*}{o[-2]}{}
\]
and the parahoric $N$ corresponds to the two vertices interchanged by Frobenius at either end of the diagram. When the dimension $n = 2m+1 \geq 5$ of $V$ is odd, we have the diagram\[
\Dynkin{A*}{e[2]}{}
\]
When the determinant of $V$ is not a norm, the parahoric $N$ corresponds to the two vertices interchanged by Frobenius at the right of the diagram.

\section{The fibers of the map $f: Y \rightarrow Z$}

In this section, we continue to assume that $k$ is a non-Archimedian local field. If $V$ is an orthogonal space over $k$ of dimension $n \geq 3$, with unit determinant and Hasse-Witt invariant $\epsilon = -1$, we have defined a parahoric subgroup $N$ up to conjugacy in $G = \SO(V)$, as well as an unramified compact subgroup $J$ up to conjugacy in $N$. Similarly to a Hermitian space $V$ over the unramified quadratic extension $k$ of $K$, whose Hermitian determinant is not a norm, we have defined a parahoric subgroup $N$ up to conjugacy in $G = U(V)$, as well as an unramified compact subgroup $J$ up to conjugacy in $N$. The groups $J$ and $N$ have the same reductive quotients. We let $Y = G/J$ and $Z = G/N$ be the associated homogeneous spaces and let $f: Y \rightarrow Z$ be the covering map. Fix a point $z \in Z$. In this section we will describe an $A_K$-analytic structure on the inverse image $f^{-1}(z) \cong N/J$ in $Y$. 

The fact that there is a natural analytic structure on the fibre, whereas $Y$ and $Z$ only have the structure of sets, follows from the fact that the compact groups $N$ and $J$ are both the $A$ valued points of smooth parahoric group schemes $\bf N$ and $\bf J$ over $A$. Since $A$ is a complete discrete valuation ring, the quotient $\bf N/\bf J$ has the structure of a formal scheme over $A$. Since $\bf J$ is a reductive group scheme over $A$, $H^1(A, \bf J) = 0$. Hence $({\bf N}/{\bf J)}(A) = {\bf N}(A)/{\bf J}(A) = N/J$. We will determine the formal scheme and analytic space structure of the fiber explicity, and see that it has an $A_K$ structure.

In the orthogonal case, fix a point $y \in Y$ with $f(y) = z$. This expresses the lattice $\Lambda$ as a direct sum $\Lambda = A_K.e + L$. To parametrize the other points in $f^{-1}(z)$, we need to determine all other such decompositions of the lattice $\Lambda$, or equivalently the rank two $A$-submodules $M = Av + Aw$ of $\Lambda$ which are isometric to $A_K.e$. The key observation is that $\pi \Lambda^{\vee}/\pi \Lambda$ is the radical of the quadratic space $\Lambda/\pi \Lambda$ over $\F$, and that a rank $2$ submodule $M$ isometric to $A_K.e$ must have this reduction modulo $\pi$. Consequently $M$ has a basis of the form
$$v = \alpha.e + \lambda$$
$$w = \beta.e + \mu$$
where $\alpha$ and $\beta$ are units of $A_K$ which are independent in the $A/\pi A$ vector space $A_K/\pi A_K$ and $\lambda$ and $\mu$ are elements of the sublattice $\pi L$. We then find inner products
$$\langle v,v \rangle = 2N(\alpha)\pi + O(\pi^2)~~~~~~\langle w,w \rangle = 2N(\beta)\pi + O(\pi^2) ~~~~~~~~~ \langle v,w \rangle = Tr(\alpha.\overline{\beta}) \pi + O(\pi^2)$$
and from this it follows that $M = Av + Aw$ is a lattice (oriented, and isometric to $A_K.e$) in a two dimensional orthogonal space over $k$ with determinant $-D$ and Hasse-Witt invariant $-1$. Taking the orthogonal complement $M^{\perp}$ of $M$ in $\Lambda$, we obtain another point $y^*$ in $f^{-1}(z)$. We still need to check that the lattice $M^{\perp}$ is non-degenerate. This argument is a bit simpler in the Hermitian case, and we postpone it for the moment.

We can parametrize the data determining the decomposition $\Lambda = M + M^{\perp}$ by the element
$$\alpha^{-1} \otimes \lambda + \beta^{-1} \otimes \mu \in A_K \otimes \pi L$$
Picking a basis of $L$ over $A$, this is just the set of points $(\pi A_K)^{n-2}$ in the formal scheme $Spf A_K[[z_1,z_2,\ldots, z_{n-2}]]$. So we have established the following.

\begin{theorem}

Assume that $V$ is an orthogonal space over $k$ of dimension $n \geq 3$ over $k$ with unit determinant and Hasse-Witt invariant $-1$. Let $Y$ be the homogeneous space parametrizing pairs $(W,L)$ consisting of an oriented two plane $W$ in $V$ with determinant $-D$ and Hasse-Witt invariant $-1$, together with a non-degenerate lattice $L$ in $W^{\perp}$. Let $Z$ be the homogeneous space of maximal lattices $\Lambda$ in $V$ together with an orientation of the lattice $\Lambda^{\vee}/\Lambda$ over the residue field of $k$.

Let $f: Y \rightarrow Z$ be the $\SO(V)$-equivariant map taking the pair $y = (W,L)$ to the non-degenerate lattice $z = \Lambda = A_K.e + L$ with orientation coming from the orientation of $W$.  Then the inverse image $f^{-1}(z)$ in $Y$ can be parametrized by the points of  $A_K \otimes \pi L \cong (\pi A_K)^{n-2}$.
\end{theorem}

We have a filtration of the parahoric $N$ by normal subgroups $N_m$ of finite index which act trivially on the quotient lattice $\Lambda/\pi^m \Lambda$. This gives a filtration of the quotient set
$$N/J > N_1J/J > N_2J/J > \ldots$$
If we identify $f^{-1}(z) = N/J$ with the points of $(\pi_K)^{n-2}$, then the subset $N_mJ/J$ corresponds to the points in $(\pi_K^{m+1})^{n-2}$. Indeed, the quotient groups
$N_mJ/N_{m+1}J$ are isomorphic to $(A_K/\pi A_K)^{n-2}$ for all $m$. 

In the simplest case, when $n = 3$ and $d \equiv 1$ the space $V$ is given by the elements of trace zero in the quaternion division algebra $B$ over $k$. Let $R$ be the unique maximal order in $B$. Then $A_K$ embeds as a subring of $R$, and we have the orthogonal decomposition
$$R = A_K + A_K.\Pi$$
where $\Pi$ is a uniformizing parameter of $R$, which normalizes $A_K$ and whose square is a uniformizing element of $A$.
The special orthogonal group $\SO(V) = B^*/k^*$ is compact, and contains the parahoric subgroup $N = R^*/A^*$ with index $2$. Hence the space $Z$ consists of two points. The subgroup $H = K^*/k^* = A_K^*/A^*$ is compact in this special case, so $H = J$ and $X = Y= B^*/A_K^*$ is isomorphic to two copies of the maximal ideal $\pi A_K$.

In the Hermitian case, a point $y \in Y$ corresponds to a line $W$ in the Hermitian space $V$ and a non-degenerate lattice $L$ in the orthogonal complement. The direct sum $\Lambda = A_K.e + L$ is a maximal Hermitian lattice in $V$ with $\Lambda^{\vee}/\Lambda = \pi^{-1}A_K.e/A_K.e$, and the point $z = f(y)$ in $Z$ remembers only the lattice $\Lambda$, not the decomposition. To find all points in $f^{-1}(z)$, we need o determine the rank $1$ sublattices $M = A_K.v$ in $\Lambda$ which are isometric to $A_K.e$.  Since $M/ \pi M = A_K.e/ \pi A_K.e$ is the radical of $\Lambda/\pi \Lambda$, the vector $v$ of $\Lambda$ must have the form
$$ v = \alpha.e + \lambda$$
where $\alpha$ is a unit in $A_K$ and $\lambda$ lies the the sublattice $\pi.L$. Let $M = A_K.v$ be the one dimensional lattice spanned by $v$, and let $M^{\perp}$ denote its orthogonal complement in the lattice $\Lambda$. We need to check that the lattice $M^{\perp}$ of rank $n-1$ is non-degenerate.

Since $\alpha$ is a unit, there is no loss of generality in assuming that the basis for $M$ has the form $v = e + \mu$ with $\mu = \pi \nu\in \pi.L$. Let $\{f_1,f_2,\ldots,f_{n-1}\}$ be an orthonormal basis for $L$ over $A_K$, and write $\nu = \sum \beta_i.f_i$. Then the vectors $-\overline{\beta_i}.e + f_i$ in $\Lambda$ give a basis for $M^{\perp}$ over $A_K$, and the Gram matrix of their inner products has determinant $\equiv 1$ modulo $\pi$. This establishes the non-degeneracy of the lattice $M^{\perp}$ (a slightly more complicated version of this argument also works in the orthogonal case).

We can parametrize the decomposition $\Lambda = M + M^{\perp}$ by the element
$$\alpha^{-1} \otimes \lambda \in \pi L.$$
Picking a basis for $L$ over $A_K$, this is the set $(\pi_K)^{n-1}$. In summary, we have shown the following.

\begin{theorem}
Assume that $V$ is a Hermitian space of dimension $n \geq 2$ over $K$ whose Hermitian determinant is not a norm. Let $Y$ be the homogeneous space parametrizing pairs $(W,L)$ consisting of a line $W$ in $V$ whose determinant is not a norm together with a non-degenerate lattice $L$ in $W^{\perp}$. Let $Z$ be the homogenous space parametrizing the non-degenerate lattices $\Lambda$ in $V$.

Let $f: Y \rightarrow Z$ be the $U(V)$-equivariant map which takes the pair $y = (W,L)$ to the lattice $z = \Lambda = A_K.e + L$. Then the inverse image $f^{-1}(z)$ in $Y$ can be parametrized by the points of $ \pi L \cong (\pi A_K)^{n-1}$.
\end{theorem}

Again, we have a filtration of the parahoric $N$ by normal subgroups $N_m$ of finite index which act trivially on the quotient lattice $\Lambda/\pi^m \Lambda$. This gives a filtration of the quotient set
$$N/J > N_1J/J > N_2J/J > \ldots$$
If we identify $f^{-1}(z) = N/J$ with the points of $(\pi_K)^{n-1}$, then the subset $N_mJ/J$ corresponds to the points in $(\pi_K^{m+1})^{n-1}$. Indeed, the quotient groups
$N_mJ/N_{m+1}J$ are isomorphic to $(A_K/\pi A_K)^{n-1}$ for all $m$.

\section{Shimura varieties associated to incoherent definite data}

Let $k$ be a totally real number field with ring of integers $A$, and fix incoherent definite orthogonal data $\{V_v\}$ over $k$ of dimension $n \geq 3$ and determinant $d$. For each place $v$ of $k$ let $V(v)$ be the neighboring global orthogonal space, and let $G(v)$ be the special orthogonal group of $V(v)$ over $k$. Let $\A$ be the ring of adeles of $k$ let $\A^f$ be the subring of finite adeles, and for a finite place $\frak p$ of $k$ let $\A^{f,\frak p}$ denote the subring of finite adeles away from $\frak p$.

To define the orthogonal Shimura variety $S$ associated to this data, we first choose a real place $v$ of $k$, and begin by defining $S$ over the complex quadratic extension $K_v \cong \bC$ of the completion $k_v = \R$. Let $L$ be an integral $A$ lattice in the neighboring orthogonal space $V(v)$, which has signature in $(n-2,2)$ at $v$ and signature $(n,0)$ at all other real places $w$ of $k$. We note that $L \otimes A_{\frak p}$ is non-degenerate for almost all primes $\frak p$. Then $M = \prod \SO(L_{\frak p})$ is an open compact subgroup of $G(v)(\A^f)$, and $M_{\frak p} = \SO(L_{\frak p})$ is a hyperspecial maximal compact subgroup of $\SO(V_{\frak p})$ for almost all primes $p$.

Let $T = \Res_{K_v/k_v} \mathbb G_m/\mathbb G_m$ be the one dimensional non-split torus over $k_v = \R$, and let $h: T \rightarrow G(v)_{\R}$ be the homomorphism described in the previous section, associated to negative oriented two dimensional subspaces of $V(v)$ over $k_v$. Then the conjugacy class $X_v$ of $h$ has the structure of two conjugate copies of the Hermitian symmetric space of $\SO(n-2,2)$. We define
$$S_M(K_v) = G(v)(k) \backslash X_v \times G(v)(\A^f)/M$$
This double coset space is the disjoint union of a finite number of connected components each of the form $\Gamma \backslash X_v$, where $\Gamma$ is an arithmetic subgroup of the algebraic group $G(v)(k)$. As such, it is a complex analytic orbifold, which has an algebraic structure (This algebraic structure is unique if $M$ is small enough, so that the subgroups $\Gamma$ have no non-trivial torsion). The complex algebraic varieties $S_M$ form a projective system, for $M' \subset M$, and one defines $S$ as the projective limit. This complex pro-scheme has a right action of the group $G(v)(\A^f)$, and $S_M$ is the quotient of $S$ by the open compact subgroup $M$.

The theory of canonical models shows that $S$ descends canonically to $k$, viewed as a subfield of $k_v = \R$ in $K_v = \bC$. We note that the descent to $k_v = \R$ is given by the anti-holomorphic involution of complex conjugation on $X_v$, and the action of the group $G(v)(\A^f)$ on $S$ is defined over $k$ on the canonical model. This defines the Shimura variety $S$ associated to the incoherent definite orthogonal data $\{V_v\}$. The irreducible components of $S$ are rational over the maximal abelian extension $E$ of $k$ with exponent $2$, and are permuted simply-transitively by the quotient $\A^*/k^*\A^{*2}$ of $G(v)(\A^f)$, which is isomorphic to the Galois group of $E$ over $k$ by the reciprocity homomorphism of global class field theory. To see that this Galois group is a quotient of $G(v)(\A^f)$, note that the spinor norm maps $G(v)(\A^f)$ onto the group $(\A^f)^*/(\A^f)^{*2}$. The quotient of this group by the group $k_+^*$ of totally positive elements in $k^*$ is the group $\A^*/k^*\A^{*2}$.

The reason that we use incoherent data to define $S$, rather than just the orthogonal group $G(v)$ with its real conjugacy class $X_v$, is that the latter depends on the choice of a real place of $k$, and only defines $S$ over a subfield of $\bC$. However, for any other real place $w$ of $k$, we have an isomorphism of adelic groups $G(w)(\A^f) \cong G(v)(\A^f)$ and an isomorphism complex varieties 
$$S_M(K_w) \cong G(w)(k) \backslash X_w \times G(w)(\A^f)/M$$
When $n = 3$, so $S_M$ is a Shimura curve, this isomorphism is due to Doi and Naganuma \cite{DN}. In the general case, it follows from results of Delgine and Borovoi on the conjugation of Shimura varieties. Therefore it is more symmetrical to use the incoherent definite orthogonal data to define $S$. We will speculate on the use of the neighboring orthogonal spaces $V(\frak p)$ for finite places of $k$ in the next section.

A similar definition works for the unitary Shimura varieties associated to incoherent definite Hermitian data $\{V_v\}$ of dimension $n \geq 1$ over the imaginary quadratic extension $K$ of the totally real field $k$. We choose a real place $v$ of $k$  and let $G(v)$ be the unitary group of the neighboring Hermitian space $V(v)$ over $k$. We begin by defining the Shimura variety $S$ over the complex quadratic extension $K_v = \bC$ of $k_v = \R$. Let $L$ be a Hermitian $A_K$ lattice in the neighboring space $V(v)$, which has signature $(n-1,1)$ at $v$ and signature $(n,0)$ at all the other infinite places $w$. Then $M = \prod U(L_{\frak p})$ is an open compact subgroup of $G(v)(\A^f)$ and $M_{\frak p} = U(L_{\frak p})$ is a hyperspecial maximal compact subgroup of $U(V_{\frak p})$ for almost all primes $\frak p$. 

Let $T = \Res_{K_v/k_v} \mathbb G_m/\mathbb G_m$ be the one dimensional non-split torus over $k_v = \R$, and let $h: T \rightarrow G(v)_{\R}$ be the homomorphism associated to a negative one dimensional subspace of $V(v)$ over $K_v$. Then the conjugacy class $X_v$ of $h$ is isomorphic to the Hermitian symmetric space of $U(n-1,1)$. We define
$$S_M(K_v) = G(v)(k) \backslash X_v \times G(v)(\A^f)/M$$
Again, this is a complex orbifold which has an algebraic structure, and the associated complex varieties $S_M$ form a projective system with limit $S$. The limit has a right action of the group $G(v)(\A^f)$, and $S_M$ is the quotient of $S$ by the open compact subgroup $M$.

The theory of canonical models shows that $S$ descends canonically to a pro-scheme over $K$, viewed as a subfield of $K_v = \bC$, and the action of the group $G(v)(\A^f)$ on $S$ is defined over $K$ on the canonical model. This is the Shimura variety $S$ associated to the incoherent definite Hermitian data $\{V_v\}$. The irreducible components of $S$ are rational over the maximal abelian extension $E$ of $K$ which is "dihedral" over $k$. They are permuted simply-transitively by the quotient $\A_K^*/(\A^*.K^* .\prod_{v|\infty}K_v^*)$ of $G(v)(\A^f)$, which is isomorphic the Galois group of $E$ over $K$ by the reciprocity homomorphism of global class field theory.

The reason that we use incoherent data to define $S$, rather than just the unitary group $G(v)$ with its real conjugacy class $X_v$, is that the latter depends on the choice of a real place of $k$, and only defines $S$ over a subfield of $\bC$. However, for any other real place $w$ of $k$, we have an isomorphism of adelic groups $G(w)(\A^f) \cong G(v)(\A^f)$ and an isomorphism complex varieties 
$$S_M(K_w) \cong G(w)(k) \backslash X_w \times G(w)(\A^f)/M$$
 This follows from general results of Deligne and Borovoi on the conjugation of Shimura varieties. Therefore it is more symmetric to use the incoherent definite Hermitian data to define $S$.

\section{The special locus modulo $\frak p$}

Let $k$ be a totally real number field, with ring of integers $A$. Let $\{V_v\}$ be definite incoherent orthogonal data of dimension $n$ and determinant $d$ for $k$. At each real place $v$ we have used the neighboring orthogonal space $V(v)$, its special orthogonal group $G(v) = \SO(V(v))$ over $k$, and the complex analytic space $X_v = G(v)(k_v)/H(v)(k_v)$ to describe the points of the Shimura variety $S_M$ over the quadratic extension $K_v \cong \bC$ of the corresponding completion $k_v = \R$.  It is reasonable to ask if something similar occurs at the finite places $\frak p$ of $k$. Namely, can we use the special orthogonal group $G(\frak p) = \SO(V(\frak p))$ of the neighboring space at $\frak p$ to say something about the $K_{\frak p}$ points of $S_M$, where $K_{\frak p}$ is the unramified quadratic extension of the completion $k_{\frak p}$? 

In this section, we will attempt to do so, under the assumptions that the residual characteristic is not equal to $2$ and that the $A$-lattice $L$ which we have used to define the open compact subgroup
$M = M_{\frak p} \times M^{\frak p}$ of $G(v)(\A^f)$ is non-degenerate at the prime $\frak p$. By this we mean that $L \otimes A_{\frak p}$ is a non-degenerate lattice in the orthogonal space $V_{\frak p}$. This assumption, which is true for almost all primes $\frak p$, implies that $d$ is a unit at $\frak p$, that $\epsilon(V_{\frak p}) = +1$, and that $M_{\frak p}$ is a hyperspecial maximal compact subgroup of $G(v)(k_{\frak p})$. In this case, it is known that the Shimura variety $S_M$, associated to the incoherent data and the choice of $M$, has a model over $A_{\frak p}$ with good reduction modulo $\frak p$ \cite{HPR}. The neighboring orthogonal space $V(\frak p)$ is positive definite at all real places of $k$. At $\frak p$ its determinant is a unit and $\epsilon_{\frak p}(V(\frak p)) = -1$. 

Let $K_{\frak p} = k_{\frak p}(\sqrt D)$ be the unramified quadratic extension of $k_{\frak p}$, let $W_{\frak p}$ be an two dimensional orthogonal space over $k_{\frak p}$ of determinant $-D$ with $\epsilon(W_{\frak p}) = -1$ and let $T$ be the one dimensional torus $(\Res_{K_{\frak p}/k_{\frak p}} \mathbb G_m)_{N=1}$. Recall that we have proved the existence and conjugacy of isometric embeddings $W_{\frak p} \rightarrow V(\frak p)$ over $k_{\frak p}$. Choosing an orientation of $W_{\frak p}$ gives a homomorphism $h: T \rightarrow G(\frak p)$ over $k_{\frak p}$ as in the real case. Writing $V(\frak p) = W_{\frak p} \oplus W_{\frak p}^{\perp}$ as an orthogonal direct sum, and choosing a non-degenerate lattice $L_{\frak p}$ in $W_{\frak p}^{\perp}$, we obtain a hyperspecial maximal compact subgroup $J_{\frak p} = \SO(W_{\frak p}) \times \SO(L_{\frak p})$ of the centralizer  of $h$. We then defined the set $Y_{\frak p}  \cong G(\frak p)(k_{\frak p})/J_{\frak p}$, which parametrizes the pairs $(W_{\frak p},L_{\frak p})$ of oriented planes of this type with the choice of a non-degenerate lattice in the orthogonal complement. 

We also defined a map $f: Y_{\frak p} \rightarrow Z_{\frak p}$ to the space of maximal lattices $\Lambda_{\frak p}$ in $V(\frak p)$, with an orientation.  Fix a vector $e$ in $W_{\frak p}$ with $q(e) = \pi$ a uniformizing element of $k_{\frak p}$ and let $A_{K_{\frak p}}.e$ be the corresponding lattice in $W_{\frak p}$. Then $\Lambda_{\frak p} = A_{K_{\frak p}}.e + L_{\frak p}$ is an oriented lattice in $V(\frak p)$, whose stabilizer is the parahoric subgroup $N_{\frak p}$ with the same reductive quotient as $J_{\frak p}$.

We would like an interpretation of the double coset spaces

$$C = G(\frak p)(k)\backslash (Z_{\frak p} \times G(\frak p)(\A^{f.\frak p}) / M^{\frak p}) \cong G(\frak p)(k)\backslash G(\frak p)(\A^f) / (N_{\frak p} \times M^{\frak p}).$$
$$D = G(\frak p)(k) \backslash (Y_{\frak p} \times G(\frak p)(\A^{f,\frak p})/M^{\frak p}) \cong G(\frak p)(k) \backslash G(\frak p)(\A^{f})/(J_{\frak p} \times M^{\frak p}).$$

Since the group $G(\frak p)(k_v)$ is compact at all the real places $v$ of $k$, and the product $(N_{\frak p} \times M^{\frak p})$ is open in $G(\frak p)(\A^{f})$, the double coset space $C$ is finite. At each point $c$ of $C$, we have a finite arithmetic subgroup $\Gamma_c$ of $G(\frak p)(k)$, which is defined by the intersection
$$\Gamma_c = G(\frak p)(k) \cap c.(N_{\frak p} \times M^{\frak p}).c^{-1}$$
The intersection takes place inside the group $G(\frak p)(\A^f)$. 

The double coset space $D$ has the structure of a $A_{K_{\frak p}}$-analytic orbifold. Via the map $F: Y_{\frak p} \rightarrow Z_{\frak p}$,we obtain a map $F_D: D \rightarrow C$. The fiber of the map $F:Y_{\frak p} \rightarrow Z_{\frak p}$ over a point $z \in Z_{\frak p}$ is isomorphic to the polydisc $(\pi A_{K_{\frak p}})^{n-2} \cong N_{\frak p}/J_{\frak p}$. Hence the fiber of the map $F_D$ over the point $c \in C$ is isomorphic to the quotient of the polydisc of dimension $n-2$ over $A_{K_{\frak p}}$ by the finite group $\Gamma_c$. There is an involution of the finite set $C$ given by the action of the normalizer of $N_{\frak p}$, and this is compatible with the involution of $D$ given by the normalizer of $J_{\frak p}$.

\begin{con}

The finite set $C$ parametrizes the set of special points on the orthogonal Shimura variety $S_M$ modulo $\frak p$. These points are all rational over the quadratic extension $\F_{\frak p^2}$ of the residue field $\F_{\frak p}$ of $k_{\frak p}$, and the action of the Galois group is given by the normalizer of $N_{\frak p}$.

The $A_{K_{\frak p}}$-analytic orbifold $D$ parametrizes the set of points of the orthogonal Shimura variety $S_M$ over $A_{K_{\frak p}}$ which have special reduction modulo $\frak p$. The reduction map is given by $F_D$ and the action of the Galois group of $A_{K_{\frak p}}/A_{\frak p}$ is given by the normalizer of $J_{\frak p}$.

\end{con}

When $n = 3$, so $S_M$ is a Shimura curve over $k$, the special points correspond to supersingular data, and the conjecture is true by results of Carayol \cite{C}.

In the Hermitian case, we expect a similar conjecture to hold, mutatis mutandis. Here we assume that $k$ is totally real, that $K$ is a totally imaginary quadratic extension of $k$, and that $\frak p$ is a finite prime of $k$ which remains inert in $K$. We start with incoherent Hermitian data, choose a real place $v$, and assume that the $A_{K}$-lattice in $V(v)$ that we have used to define the open compact subgroup $M= M_{\frak p} \times M^{\frak p}$ is non-degenerate. This implies that the Hermitian determinant $d$ is a norm locally at the prime $\frak p$ and that $M_{\frak p}$ is a hyperspecial compact subgroup of $U(v)(k_{\frak p})$. In this case, the unitary Shimura variety $S_M$, associated to the incoherent data and the choice of $M$, is defined over $K$ and has a model over $A_{K_{\frak p}}$ with good reduction modulo $p$ \cite{RSZ}. The neighboring Hermitian space $V(\frak p)$ is positive definite at all complex places of $K$. At the prime $\frak p$ its Hermitian determinant is not a norm. We would like to use the corresponding global unitary group $G(\frak p) = U(V(\frak p))$ and the two local homogeneous spaces $Z_{\frak p}$ and $Y_{\frak p}$ for $G(\frak p)(k_{\frak p})$ to parametrize the finite set of special points over $\F_{\frak p^2}$ as well as the points in $A_{K_{\frak p}}$ which reduce to special points modulo $\frak p$.

As in the orthogonal case, we consider the double coset spaces

$$C = G(\frak p)(k)\backslash (Z_{\frak p} \times G(\frak p)(\A^{f.{\frak p}}) / M^{\frak p}) \cong G(p)(k)\backslash G(p)(\A^f) / (N_p \times M^p).$$
$$D = G(\frak p)(k) \backslash (Y_{\frak p} \times G(\frak p)(\A^{f,\frak p})/M^{\frak p}) \cong G(\frak p)(k) \backslash G(\frak p)(\A^f)/(J_{\frak p} \times M^{\frak p}).$$

The double coset space $C$ is finite, and at each point $c$ of $C$, we have a finite arithmetic subgroup $\Gamma_c$ of $G(\frak p)(k)$, which is defined by the intersection
$$\Gamma_c = G(\frak p)(k) \cap c.(N_{\frak p} \times M^{\frak p}).c^{-1}$$
The intersection takes place inside the group $G(\frak p)(\A^f)$. 

This double coset space $D$ has the structure of a $A_{K_{\frak p}}$-analytic orbifold. Via the map $F: Y_{\frak p} \rightarrow Z_{\frak p}$, we obtain a map $F_D: D \rightarrow C$.
Since the fiber of the map $F:Y_{\frak p} \rightarrow Z_{\frak p}$ over a point $z \in Z_{\frak p}$ is isomorphic to the polydisc $(\pi A_{K_{\frak p}})^{n-1} \cong N_{\frak p}/J_{\frak p}$, the fiber of the map $F_D$ over the point $c \in C$ is isomorphic to the quotient of the polydisc of dimension $n-1$ over $A_{K_{\frak p}}$ by the finite group $\Gamma_c$.

\begin{con}

The finite set $C$ parametrizes the set of special points on the unitary Shimura variety $S_M$ modulo $\frak p$. These points are all rational over the residue field $\F_{\frak p^2}$ of $\frak p$ in $A_K$.

The $A_{K_{\frak p}}$-analytic orbifold $D$ parametrizes the set of points of the unitary Shimura variety $S_M$ over $A_{K_{\frak p}}$ which have special reduction modulo $\frak p$. The reduction map is given by $F_D$.

\end{con}

Finally, note that when $M_{\frak p}$ is hyperspecial, the components of the orthogonal Shimura variety $S_M$ are all rational over the quadratic extension $K_{\frak p}$ of the localization $k_{\frak p}$. This is certainly necessary for the existence of any $K_{\frak p}$ rational points, let alone points reducing to the special locus! Indeed in the orthogonal case, the components of $S_M$ are rational over a finite abelian extension $E$ of $k$ of exponent $2$. This extension depends on the image of the spinor norm from the open compact subgroup $M$ of $G(v)(\A^f)$. When $M_{\frak p}$ is hyperspecial, the image of the spinor norm contains the unit classes in $k_{\frak p}^*/k_{\frak p}^{*2}$ and the extension $E/k$ is unramified at the place $\frak p$. Hence the decomposition group of $ \frak p$ in the Galois group $\Gal(E/k)$ is cyclic, so has order $1$ or $2$, and the components of $S_M$ are all rational over the unramified quadratic extension $K_{\frak p}$ of $k_{\frak p}$. A similar result holds in the Hermitian case, as a prime $\frak p$ which is inert in the quadratic extension $K$ where $M_{\frak p}$ is hyperspecial splits completely in the ring class extension $E$ of $K$ where the components of $S_M$ are rational.

\section{A mass formula}

The sum over the double coset space $C$ (which conjecturally parametrizes the special points of $S_M$ modulo $\frak p$) of the reciprocals of the order of the corresponding finite groups
$\Gamma_c$ 
$$\Mass(C) =\sum_{c \in C} 1/\#\Gamma_c$$
is given by a "mass formula", using the $L$-function of the motive and the Tamagawa number of $G(\frak p)$ (cf. \cite{G4}, \cite{GHY}). 

Since the motive and Tamagawa number of $G(\frak p)$ are equal to the motive and Tamagawa number of its inner form $G(v)$ over $k$, there is a simple relation between this mass and the virtual Euler characteristic $\chi$ of the complex orbifold $S_M(K_v)$, as defined by Serre \cite{S2}
$$\#(W_G/W_H) \times \Mass(C)=  \#G_0(\F_{\frak p})\big/\#N_0(\F_{\frak p}).(-q)^{\dim S} \times \chi(S_M(K_v)).$$
where $W_G$ denotes the Weyl group of $G = G(v)$, $W_H$ denotes the Weyl group of the centralizer of $h$, and  $q = \#F_{\frak p}$.

For example, when $\{V_v\}$ is incoherent definite orthogonal data of odd dimension $2n+1$, we find the formula
$$2n \times \Mass(C) = (1 - q + q^2 - q^3 + \ldots + (-q)^{2n-1}) \times \chi(S_M(K_v)).$$
In the case when $n = 1$, $S_M$ is a Shimura curve. If we assume that $M$ is small enough so that all the finite groups $\Gamma_c$ are trivial, the mass is just the number of supersingular points on the curve $S_M$ modulo $\frak p$. Since these points are all rational over the quadratic extension of $\F_{\frak p}$, which is the field of $q^2$ elements, the formula implies
$$2 \times \#S_M(q^2) \geq (1-q) \times (2 - 2g)$$
where $g$ is the genus of the complex curve $S_M$. Hence the tower of curves $S_M$ realize the asymptoptic Drinfeld-Vladut bound \cite{DV} of $q-1$, for the ratio of the number of points for a curve over a field with $q^2$ elements to its genus $g$.

\def\noopsort#1{}
\providecommand{\bysame}{\leavevmode\hbox to3em{\hrulefill}\thinspace}

\end{document}